\documentclass[11pt]{amsart}

\usepackage{amscd}
\usepackage{amssymb}
\usepackage{boxedminipage}
\usepackage{enumerate}
\usepackage{euscript}

\theoremstyle{plain} \newtheorem{theorem}{Theorem}[section]
\theoremstyle{plain} \newtheorem{lemma}[theorem]{Lemma}
\theoremstyle{plain} \newtheorem{proposition}[theorem]{Proposition}

\newtheorem{corollary}[theorem]{Corollary}
\newtheorem{problem}{Problem}

\newcommand{\nr}{\refstepcounter{theorem}  
                   \noindent {\thetheorem .}}
\newcommand{\defi}{\medskip \noindent {\it Definition \nr} }
\newcommand{\defifin}{\medskip}
\newcommand{\eks}{\medskip \noindent {\it Example \nr} }
\newcommand{\eksfin}{\medskip}
\newcommand{\rem}{\medskip \noindent {\it Remark \nr} }
\newcommand{\remfin}{\medskip}

\newcommand{\note}{\medskip \noindent {\it Notation \nr} }
\newcommand{\notefin}{\medskip}

\newcommand{\llabel}{\addtocounter{theorem}{-1}
\refstepcounter{theorem} \label}

\newcommand{\op}{{\mathcal O}}

\newcommand{\go}{\op}

\newcommand{\she}{\EuScript{S}\text{h}}
\newcommand{\cm}{\EuScript{CM}}
\newcommand{\cmd}{\EuScript{CM}^\dagger}
\newcommand{\cmri}{\EuScript{CM}^\circ}
\newcommand{\cler}{\EuScript{CL}}
\newcommand{\clerd}{\EuScript{CL}^\dagger}
\newcommand{\clerri}{\EuScript{CL}^\circ}
\newcommand{\gor}{\EuScript{G}}
\newcommand{\gF}{\mathcal{F}}
\newcommand{\gM}{\mathcal{M}}
\newcommand{\gE}{\mathcal{E}}
\newcommand{\gD}{\mathcal{D}}
\newcommand{\gI}{\mathcal{I}}
\newcommand{\gP}{\mathcal{P}}
\newcommand{\gK}{\mathcal{K}}
\newcommand{\gS}{\mathcal{S}}

\newcommand{\res}{\text{res}}
\newcommand{\Hom}{\text{Hom}}

\newcommand{\Tor}{\text{Tor}}
\newcommand{\ghom}{\mathcal{H}om}
\newcommand{\gext}{\mathcal{E}xt}

\newcommand{\rank}{\text{rank}\,}

\newcommand{\sus}{\subseteq}
\newcommand{\sups}{\supseteq}
\newcommand{\pil}{\rightarrow}
\newcommand{\vpil}{\leftarrow}

\newcommand{\inpil}{\hookrightarrow}

\newcommand{\mto}[1]{\stackrel{#1}\longrightarrow}
\newcommand{\vmto}[1]{\stackrel{#1}\longleftarrow}

\newcommand{\iso}{\cong}
\newcommand{\te}{\otimes}

\newcommand{\tL}{\tilde{L}}
\newcommand{\tM}{\tilde{M}}

\newcommand{\tvH}{\widetilde{H}}
\newcommand{\tvh}{\widetilde{h}}

\newcommand{\tS}{\tilde{S}}
\newcommand{\tT}{\tilde{T}}
\newcommand{\tR}{\tilde{R}}

\newcommand{\ts}{\tilde{s}}

\newcommand{\empt}{\emptyset}
\newcommand{\bfa}{{\bf a}}

\newcommand{\dl}{\Delta}

\newcommand{\lk}{\text{lk}}
\newcommand{\lkd}{\lk_\Delta}

\newcommand{\pnm}{{\bf P}^{n-1}}
\newcommand{\opnm}{{\go_{\pnm}}}
\newcommand{\ompnm}{\omega_{\pnm}}

\newcommand{\hele}{{\bf Z}}

\def\PP{{\mathbb P}}

\begin{document}
\title [Hierarchies of simplicial complexes]
{Hierarchies of simplicial complexes via the BGG-correspondence}
\author { Gunnar Fl{\o}ystad}
\address{ Matematisk Institutt\\
          Johs. Brunsgt. 12 \\
          5008 Bergen \\
          Norway}   
        
\email{ gunnar@mi.uib.no}

\begin{abstract}
Via the BGG-correspondence a simplicial complex $\dl$ on $[n]$ is 
transformed into a complex of coherent sheaves $\tL(\dl)$ on the projective
space $\pnm$. In general we compute the support of each of its cohomology 
sheaves.

When the Alexander dual $\dl^*$ is Cohen-Macaulay there is only one
such non-zero cohomology sheaf. We investigate when this sheaf can be
an $a$'th syzygy sheaf in a locally free resolution
and show that this corresponds exactly to the case of $\dl^*$ being 
$a+1$-Cohen-Macaulay as defined by K.Baclawski \cite{Ba}.


By putting further conditions on the sheaves we get nice
subclasses of $a+1$-Cohen-Macaulay simplicial complexes whose 
$f$-vector depends only on $a$ and the invariants
$n,d$, and $c$. When $a=0$ these are the bi-Cohen-Macaulay simplicial
complexes, when $a=1$ and $d=2c$ cyclic polytopes are examples, 
and when $a=c$ we get 
Alexander duals of the Steiner systems $S(c,d,n)$.

We also show that $\dl^*$ is Gorenstein* iff the associated coherent sheaf
of $\dl$ is an ideal sheaf.

\end{abstract}

\maketitle

\section*{Introduction}

A simplicial complex $\dl$ on $\{1, \ldots, n\}$ corresponds to a monomial 
ideal in the exterior algebra. Via the BGG-correspondence such an ideal 
transforms to a complex of coherent sheaves $\tL(\dl)$ on the projective
space $\pnm$. In a recent paper, \cite{FV}, J.E.Vatne and the author studied
when the BGG-correspondence applied to a simplicial complex $\dl$ gives rise
to a complex $\tL(\dl)$ with at most one non-zero cohomology sheaf
$\gS(\dl)$. We showed that this happens exactly when the Alexander dual 
$\dl^*$ is a Cohen-Macaulay simplicial complex. We then singled out a 
particularly nice class of coherent sheaves, the locally Cohen-Macaulay
sheaves and asked when $\gS(\dl)$ belonged to this class. This turned
out to happen when both $\dl$ and $\dl^*$ are Cohen-Macaulay.

In this paper we study this further in several directions. Firstly we develop
in Section 3 criteria for telling which cohomology sheaves of $\tL(\dl)$
are non-zero and if so, what their support is.
We also describe the multigraded Hilbert functions of
these cohomology sheaves. This is related to the cellular complexes
of \cite{BS}.

Secondly, we again single out nice classes of coherent sheaves and ask
when a simplicial complex $\dl$ is such that $\gS(\dl)$ is in such a class.
This time we consider the class of coherent sheaves $\she_a$ which can
occur as $a$'th syzygy sheaf in a locally free resolution of a coherent 
sheaf. Thus when $a$ is $1$ we have the torsion free sheaves, when $a$ is $2$ 
the reflexive sheaves and when $a$ is $n-1$ the vector bundles. In Section 4
we consider the class $\cler_a$ of simplicial complexes $\dl$ such that
$\gS(\dl)$ is in $\she_a$, giving us a "hierarchy" 
$\cler = \cler_0 \sups \cler_1 
\sups \cdots$ of simplicial complexes. Note that $\cler$ is the class of 
Alexander duals of Cohen-Macaulay simplicial complexes. There turns out
to be nice, compact criteria for membership in the classes $\cler_a$ 
and when taking restrictions and links of a $\dl$ in $\cler_a$ these
are predictably in $\cler_a$ and $\cler_{a-1}$ respectively. 
(Via Alexander duals this 
generalizes the fact that the link of a Cohen-Macaulay simplicial complex
is Cohen-Macaulay.) 

The three most important invariants of a simplicial complex $\dl$ is $n$, its 
dimension $d-1$ and $c$,
the largest integer such
that all $c$-sets are in $\dl$. Via certain nice additional properties 
which one may put on sheaves in $\she_a$ we further single out a particular
nice subclass of $\cler_a$. When $a$ is $0$ this subclass
turns out to be the class
studied in \cite{FV}, the bi-Cohen-Macaulay simplicial complexes. When
$a$ is $c$, the maximal interesting value for $a$, it turns out that we 
get exactly the Steiner systems $S(c,d,n)$. In general for each $a$, the
simplicial complexes of this nice subclass of $\cler_a$ will have $f$-vectors
depending only on $c,d$ and $n$.

\medskip

Now taking Alexander duals of the classes $\cler_a$ we get in Section 5
a "hierarchy" of 
Cohen-Macaulay simplicial complexes $\cm = \cm_0 \sups \cm_1 \sups \cdots$.
The simplicial complexes in $\cm_a$ turn out to be exactly the 
$a+1$-Cohen-Macaulay simplicial complexes as defined by K. Baclawski \cite{Ba}.
(Thus the simplicial complexes in $\cler_a$ are the Alexander duals of the
$a+1$-Cohen-Macaulay simplicial complexes.)
After translating most of the results for the
classes $\cler_a$ and its subclasses to the classes $\cm_a$ 
and its subclasses, we consider the subclass of $\cm$
consisting of Gorenstein* simplicial complexes $\dl$. This 
class corresponds exactly to the subclass of $\cm_1$ with 
$\tvH_{d-1}(\dl) = k$ and 
is the class of Cohen-Macaulay simplicial complexes $\dl$
such that $\gS(\dl^*)$ is a torsion free rank one coherent sheaf. Moreover
we show that this sheaf may be naturally identified with the 
associated coherent sheaf of 
the ideal defining the Stanley-Reisner ring of $\dl^*$.

Having described the contents of Sections 3,4, and 5 we inform that 
Section 1 contains preliminaries and techniques for dealing with simplicial 
complexes. Section 2 describes the BGG-correspondence and the classes
of coherent sheaves we single out for study, namely the classes of sheaves
which can be $a$'th syzygy sheaves in a locally free resolution.
In the last section, Section 6, we pose some problems for further study.

\section{Simplicial complexes.}

Denote by $[n]$ the set of integers $\{1, \ldots, n\}$. A simplicial complex 
$\dl$ on $[n]$ is a family of subsets of $[n]$ such that if $Y \sus X \sus [n]$
and $X$ is in $\dl$ then $Y$ is in $\dl$. 

\subsection{Notions}
We recall some notions for simplicial complexes.
An element in $\dl$ is called
a {\it face} of $\dl$, and a maximal face is a {\it facet}.
If $d$ is the maximal cardinality of a face
of $\dl$, the {\it dimension} of $\dl$ is $d-1$. If $c$ is the largest integer
such that all $c$-sets in $[n]$ are in $\dl$, we call
$c-1$ the {\it frame dimension} of $\dl$. By convention the empty simplex
$\emptyset$ has $c = -1$ (while $\{\emptyset\}$ has $c=0$).

The {\it Alexander dual} simplicial complex $\dl^*$ is the simplicial complex
on $[n]$ consisting of all $F$ in $[n]$ such that the complement
$F^c = [n]\backslash F$ is not in $\dl$.
Let $d^*-1$ and $c^*-1$ be the dimension and frame
dimension of $\dl^*$. It is easily seen that
\[ n = d + c^* + 1, \quad n = d^* + c + 1.\]

Now we introduce the following notation. Let $R \cup S \cup T$
be a partition of $[n]$ into three disjoint subsets, and let
\[ \dl^{S,T}_R = \{F \sus R \, | \,  F \cup S \in \dl\} \]
which is a simplicial complex on $R$. Note that 
$\dl^{S, \empt}_{S^c}$ is the link $\lkd S$, and
that $\dl^{\empt, R^c}_R$ is the restriction $\dl_R$ to $R$.
Also note that $\dl^{S,T}_R$ is  $\empt$ iff $S$ is not a face of $\dl$.

\begin{lemma} \label{1alexdual}

\[(\dl^{S,T}_R)^* = (\dl^*)^{T,S}_R.\]

\end{lemma}

\begin{proof}
Let $F$ be a subset of $R$. That $F$ is in $(\dl^{S,T}_R)^*$ means that
$R\backslash F$ is not in $\dl^{S,T}_R$ or $(R\backslash F) \cup S$ not in 
$\dl$.

That $F$ is in $(\dl^*)^{T,S}_R$ means that $F \cup T$ is in $\dl^*$ or
$[n]\backslash (F \cup T)$ not in $\dl$. But the latter is 
$(R\backslash F) \cup S$.
\end{proof}

\subsection{Homology}
Let $V$ be a vector space over a field $k$ with basis $e_1, \ldots, e_n$ 
and $E= E(V)$ the exterior
algebra $\oplus_0^n \wedge^i V$. Let $W = V^*$ be the dual space with dual 
basis $x_1, \ldots, x_n$.
Consider the monomials $e_{i_1}\cdots e_{i_r}$ in $E$ such that
$\{i_1, \ldots, i_r\}$ is not in $\dl$. They form a basis for an ideal
$J_\dl$ in $E$.


Dualizing the inclusion $J_\dl \sus E$ we get an exact sequence of
$E$-modules
\begin{equation} 
0 \pil C_\dl \pil E(W) \pil (J_\dl)^* \pil 0. \label{1cdel}
\end{equation}
where $C_\dl$ is the kernel.
A basis for $C_\dl$ consist of all monomials $x_{i_1} \cdots x_{i_r}$
such that $\{i_1, \ldots, i_r\}$ are in $\dl$.
Left multiplication with $u = e_1 + e_2 + \cdots + e_n$ gives a differential
$d$ on $C_\dl$ and the reduced homology of $\dl$ is defined by 
($W$ has degree $-1$)
\[ \tvH_p(\dl) = H^{-p-1} (C_\dl, d). \]

Note that via the isomorphism of $E$-modules between 
$E(V)$ and $E(W)(-n)$, the ideal $J_\dl$ gets identified with the 
submodule $C_{\dl^*}(-n)$ of $E(W)(-n)$. We therefore get an exact sequence
\begin{equation}
 0 \pil C_\dl \pil E(W) \pil (C_{\dl^*})^*(n) \pil 0. \label{1check}
\end{equation}

\begin{lemma} \label{1alexdualhom}
$\tvH_{n-3-p}(\dl^*) \iso \tvH_p(\dl)^*.$
\end{lemma}

\begin{proof}
Since the differential on $E(W)$ is acyclic we get from (\ref{1check})
\[ H^{-p-1}(C_\dl) \iso H^{n-2-p}((C_{\dl^*})^*) 
= H^{p+2-n}(C_{\dl^*})^*. \]
\end{proof}

\subsection{Reducing}

Now if $[n]$ is partioned as $R \cup \{x\}$ we get a short exact sequence
of $E$-modules
\[ 0 \pil C_{\dl_R} \pil C_{\dl} \pil C_{\lk_\dl \{x\}}(1) \pil 0. \]
The following is a basic kind of general deduction from this type
of exact sequence. 

\begin{proposition}
Let $R \cup S \cup T$ be a partition of $[n]$ and suppose 
$\tvH_p(\dl_R^{S,T})$ is non-zero.

a) (Reducing $S$.) Given $S^\prime \sus S$. Then there are $R^\prime \sups R$ 
and $T^\prime \sups T$ such that 
$\tvH_{p^\prime}(\dl_{R^\prime}^{S^\prime, T^\prime})$ is non-zero,
where $p^\prime - p$ is the cardinality of $R^\prime \backslash R$.

b) (Reducing $T$.) Given $T^\prime \sus T$. Then there are $R^\prime \sups R$ 
and $S^\prime \sups S$ such that 
$\tvH_{p}(\dl_{R^\prime}^{S^\prime, T^\prime})$ is non-zero.

c) (Reducing $R$.) Given $R^\prime \sus R$. Then there are $S^\prime \sups S$ 
and $T^\prime \sups T$ such that 
$\tvH_{p^\prime}(\dl_{R^\prime}^{S^\prime, T^\prime})$ is non-zero
where $p - p^\prime$ is the cardinality of  $S^\prime \backslash S$.

\end{proposition}

\begin{proof} All of these just follow from the exact sequence
\[ 0 \pil C_{\dl_R^{S, T\cup\{x\}}} \pil C_{\dl_{R \cup \{x\}}^{S,T}}
\pil C_{\dl_R^{S \cup \{x\},T}}(1) \pil 0 \]
where $R \cup S \cup T \cup \{x\}$ is a partition of $[n]$, by running long
exact cohomology sequences.
\end{proof}

\rem This proposition generalizes Corollary 4.4 of \cite{AAH}.
\remfin

\note We shall usually write $\dl^{S,T}_{-R}$ for 
$\dl^{S,T}_{[n]\backslash R}$ (in which case $R = S \cup T$). Often we shall
also drop the $T$ and write this simply as $\dl^S_{-R}$.
In general when $Y$ is a set we shall use the lower case letter $y$
to denote the cardinality of $Y$. For instance given $\dl^S_{-R}$ then
$r$ and $s$ will be the cardinalities of $R$ and $S$ respectively.
\notefin

\begin{corollary}  \label{1redcor}
Suppose $\tvH_p(\dl^{S,T}_{-S\cup T})$ is non-zero. 
Then $c-1 \leq p+s \leq d-1$ and
the following holds.

a) If $p+s = c-1$  and $S^\prime \sus S$ then 
$\tvH_{p^\prime}(\dl^{S^\prime,T}_{-S^\prime\cup T})$ is non-zero,
where $p^\prime + s^\prime = c-1$.

b) If $p+s = d-1$ and $T^\prime \sus T$ then 
$\tvH_p(\dl^{S,T^\prime}_{-S\cup T^\prime})$ is non-zero.

c) If $p \geq 0$ then $\tvH_{p-1}(\lkd S^\prime)$ is non-zero for some 
$S^\prime$ strictly containing $S$.


d) $\tvH_{p^\prime}(\dl_R)$ is non-zero for some $R$ and $p^\prime \geq p$.

e) $\tvH_p(\lkd S^\prime)$ is non-zero for some $S^\prime$ containing $S$.

\end{corollary}

\begin{proof} The frame dimension of 
$\dl^{S,T}_{-S\cup T}$ is $\geq c-s-1$. 
If therefore
$p < c-s-1$ then $\tvH_p(\dl^{S,T}_{-S\cup T})$ is zero. Similarly the 
dimension of 
$\lkd S$ is $\leq d-s-1$ and so if $p \geq d-s$ then 
$\tvH_p(\dl^{S,T}_{-S\cup T})$ is zero.

a) Reduce $S$ to $S^\prime$ and get $\tvH_{p^\prime}
(\dl^{S^\prime,T^\prime}_{-S^\prime\cup T^\prime})$ non-zero. First note
we must have $p^\prime + s^\prime \geq c-1$. Also we have $p^\prime -p = 
s+t - s^\prime - t^\prime$. Since $t^\prime \geq t$ we must have 
$T = T^\prime$ and $p^\prime + s^\prime =  c-1$.
Part b) is shown similiarly as a).

c) We have $\tvH_p(\dl^{S,T}_{R})$ non-zero. Now reduce $R$. 
Sooner or later
$S$ increases to $\tS = S \cup \{x\}$ such that 
\begin{eqnarray} \tvH_{p-1}(\dl^{\tS,\tT}_{\tR}) \neq 0. \label{1cp}
\end{eqnarray}
This is so since otherwise we would get 
$\tvH_p( \dl^{S, S^c}_{\emptyset})$ non-zero which is impossible. 
Now reducing $\tT$ in (\ref{1cp}) to $\emptyset$ we 
get $\tvH_{p-1}( \dl^{S^\prime, \emptyset}_{-S^{\prime}})$ non-zero.


d) and e) follow by reducing $S$, respectively $T$, to the empty set.
\end{proof}

For later we need the following.

\begin{lemma} If $\dl$ is not the $n-1$-simplex or the empty simplex, 
then $H_{c-1}(\lkd S)$ is non-zero for some $S$. \label{1lem-clk}
\end{lemma}

\begin{proof} Let $R$ be a $c+1$-set which is not a face of $\dl$.
Then clearly $H_{c-1}(\dl_R)$ is non-zero. Now reduce $T (= R^c)$ to 
$\emptyset$, and get $H_{c-1}(\lkd S)$ non-zero.
\end{proof}

\subsection {CLeray and Cohen-Macaulay simplicial complexes}

For an integer $e$, a simplicial complex $\dl$ is called $e$-Leray,
\cite[p.12]{Ka}, if $\tvH_p(\dl_R) = 0$ for all $p \geq e$ and subsets $R$
of $[n]$. Motivated by this we make the following definition.

\defi \llabel{1cler} A simplicial complex $\dl$ is {\it CLeray} if
$\tvH_p(\dl_R) = 0$ for all $p \geq c$, that is for all $p$ greater than
the frame dimension of $\dl$, and subsets $R$ of $[n]$.
(Note that we do not speak of a complex being say, 3Leray or 2Leray. One 
can speak of a complex being CLeray or 3-Leray.)
\defifin 

Note that when $c=0$ this gives that $\dl$ is a simplex on the vertices
it contains.

Recall that $\dl$ is {\it Cohen-Macaulay} 
if $\tvH_p(\lkd S)= 0$ when $p+s \leq d-2$.
The class of CLeray and Cohen-Macaulay simplicial complexes are now seen
to be Alexander dual.

\begin{proposition} $\dl$ is CLeray iff $\dl^*$ is Cohen-Macaulay.
\end{proposition}

\begin{proof} By Lemmata \ref{1alexdual} and \ref{1alexdualhom}
\[ \tvH_p(\dl_R)^* = \tvH_{r-3-p}((\dl_R^{\empt, R^c})^*) 
= \tvH_{r-3-p}((\dl^*)^{R^c, \empt}_R) = \tvH_{r-3-p}(\lk_{\dl^*} R^c).
\]
The condition that $p \geq c$ is equivalent to
$n-3-p \leq d^* -2$ and so we get the statement.
\end{proof}

Another description of $\dl$ being CLeray is the following.

\begin{proposition} \label{1leraybes} $\dl$ is CLeray iff 
$\tvH_c(\lkd S) = 0$ for all $S$ in $\dl$.
\end{proposition}

\begin{proof}
Suppose $\dl$ is not CLeray. Then $\tvH_p(\dl_R)$ is non-zero for some
$R$ and $p \geq c$. By Corollary \ref{1redcor} e), 
we get $\tvH_p(\lkd S)$ non-zero for some $S$ and 
so $\tvH_c(\lkd S^\prime)$ non-zero for some $S^\prime$ in $\dl$ by 
\ref{1redcor} c).

Conversely, if $\tvH_c(\lkd S)$ is non-zero, then by Corollary \ref{1redcor} 
d), $\tvH_p(\dl_R)$ is non-zero for some $p \geq c$ and so $\dl$ is not
CLeray.
\end{proof}

\begin{corollary}
$\dl$ is Cohen-Macaulay iff $\tvH_p(\dl_{-R})$ is zero for 
$p+r = d-2$.
\end{corollary}

\begin{proof} 
Using Lemma \ref{1alexdualhom} this follows from the above.
\end{proof}

\section{The BGG-correspondence and a hierarchiy of coherent sheaves}

Recall that $W = V^*$ is the dual space of $V$ 
and let $S = S(W)$ be the symemtric algebra on $W$. 
If $M = \oplus_{i \in \hele} M_i$ is a graded module over $E$ we can form the 
complex
\[ L(M) \, : \, \cdots \pil S(i) \te_k M_i \mto{d^i} S(i+1) \te_k 
M_{i+1} \pil \cdots  \]
with differential given by
\[ d^i (s \te m) = \sum_{\alpha = 1}^n sx_\alpha \te e_\alpha m. \]

Sheafifying this we get a complex of coherent sheaves
\[ \tL(M) \, : \, \cdots \pil \opnm(i) \te_k M_i \mto{d^i} 
\opnm(i+1) \te_k M_{i+1} \pil \cdots  \]
on the projective space $\pnm$. This, in short, is the BGG-correspondence 
between graded
modules over the exterior algebra and complexes of coherent sheaves on $\pnm$,
originally from \cite{BGG}. Our main reference for this will be \cite{EFS}.

Let us give som properties of this correspondence.

\subsection{Restriction to linear subspaces} \label{2res}

Let $V^\prime$ be a subspace of $V$ so $W^\prime = V^{\prime *}$ is a quotient 
space of $W$.
Via $E(V^\prime) \pil E(V)$, the module $M$ may be restricted to a module 
$\res\, M$ over
$E(V^\prime)$. Also via $S(W) \pil S(W^\prime)$ we may form the quotient 
complex
$L(M) \te_S S(W^\prime)$. These are related by
\[ L^\prime(\res M) = L(M) \te_S S(W^\prime) \]
where $L^\prime$ is the corresponding functor for $E(V^\prime)$-modules.
\subsection{Duals of complexes}

Now we consider $\wedge^n W$ to be a module over $E$ in degree $-n$. 

Let $M^{\vee} = \Hom_k(M,\wedge^n W).$ Then we have canonically
\[ L(M^{\vee}) = \Hom_S(L(M), S(-n) \te_k \wedge^n W)[n]. \]
Sheafifying this, note that the canonical sheaf $\ompnm$ naturally 
identifies
with the sheafification of $S(-n) \te_k \wedge^n W$, we get
\begin{equation}
 \tL(M^\vee) = \ghom_{\opnm}(\tL(M), \ompnm)[n]. \label{2dualkx}
\end{equation}
In particular if the only nonzero cohomology of $\tL(M)$ is 
$\gF = H^{-a} \tL(M)$,
a coherent sheaf, we see that
\begin{equation} 
\gext^r_{\opnm}(\gF, \omega_{\pnm}) = H^{a+r-n} \tL(M^\vee). 
\label{2extdual}
\end{equation}

\subsection {Calculating Tor's} \label{2Tor}
For a graded module $M$ over $E(V)$ and $v$ an element of $V$, let
$H^p(M,v)$ be the cohomology of the complex
\[  M_{p-1} \mto{\cdot v} M_p \mto{\cdot v} M_{p+1}. \]
Letting $k(v)$ be $\opnm / \gI_v$ where $\gI_v$ is the ideal sheaf of 
the point $v$ in $\pnm$ 
we have, due to \cite{EPY},
\[ H^p(M,v) = H^p(\tL(M) \te_{\opnm} k(v)). \]
In particular, if $\gF = H^0 \tL(M)$ is the only non-zero cohomology
\[ H^{-p}(M,v) = \Tor_p^{\opnm}(\gF, k(v)). \]

 This is what will be needed about the BGG-correspondence. Now we turn to 
describing some
classes of coherent sheaves.
 
\subsection {A hierarchy  of coherent sheaves} \hfill


\defi  A coherent sheaf $\gF$ on $\pnm$ is an {\it $a$'th syzygy sheaf} if 
there is a 
locally free resolution of some coherent sheaf $\gM$
\begin{eqnarray}
 \gM \vpil \gE_0  \vpil \ldots \vpil \gE_{i-1} \vmto{d_{i}} 
\gE_{i} \vpil \ldots \label{2resM} \end{eqnarray}
such that $\gF$ is the kernel of $d_{a-1}$. 
\defifin

Now let $\she = \she^0$ be the category of all coherent sheaves on $\pnm$.
For $a \geq 1$ let $\she_a$ be the full subcategory of $\she$ consisting 
of the sheaves which can occur as an $a$'th syzygy sheaves. There is then a 
filtration into ``nicer and nicer'' classes of coherent sheaves
\[ \she_0 \sups \she_1 \sups \cdots \sups \she_{n-1}. \]
Note that $\she_1$ are the torsion free sheaves, $\she_2$ are the
reflexive sheaves, i.e. sheaves $\gF$ such that the natural map
$\gF \pil \gF^{\vee \vee}$ is an isomorphism (here
$\gF^\vee$ is $\ghom_\opnm(\gF, \ompnm)$), and $\she_{n-1}$ are the 
vector bundles.

\begin{proposition}
$\gF$ is in $\she_a$ if and only if the codimension of the support of
$\gext^i_\opnm(\gF, \ompnm)$ is greater or equal to $i+a$ for all $i > 0$.
\end{proposition}
 
\begin{proof} If $\gF$ is in $\she_a$, then with the notation in (\ref{2resM}),
$\gext^i(\gF, \ompnm)$ is equal to $\gext^{i+a}(\gM, \ompnm)$ for $i > 0$
and the statement follows by \cite[20.9]{Ei}.

For later use note the following.

\medskip

\noindent{\bf Fact.} 
If $\gF$ is a coherent sheaf with support in 
codimenison $r$, then the sheaf $\gext^i(\gF, \ompnm)$ has codimension $r$ 
when $i=r$ and is zero when $0 < i < r$.
This follows for instance from \cite[18.4]{Ei} and its proof.
\medskip

Now assume $\gext^i(\gF, \ompnm)$ has codimension greater or equal
to $i+a$ for $i \geq 1$.
If $a \geq 1$ it follows from the above fact that $\gF$ cannot have a torsion 
subsheaf, and so $\gF$ is in $\she_1$.

Now let $a \geq 2$. Since $\gF$ is in $\she_1$, there is an exact sequence
\[ 0 \pil \gF \pil \gF^{\vee \vee} \pil \gP \pil 0\]
where $\gP $ is the cokernel. Let $1+r$ be the codimension of $\gP$. Since
any torsion free sheaf is locally free outside a codimension $2$ subset,
 $1+ r$ is $\geq 2$. Now we get a long exact sequence
\begin{eqnarray*} && \gext^r(\gF^{\vee \vee}, \ompnm) \pil \gext^r(\gF, \ompnm)
\pil \gext^{r+1}(\gP, \ompnm) \\
&\pil& \gext^{r+1}(\gF^{\vee \vee}, \ompnm).
\end{eqnarray*}
Hence if $\gP$ is non-zero, that is $\gF$ is not in $\she_2$, 
then $\gext^r(\gF, \ompnm)$ has codimension $1+r$, 
which is against our assumption. Hence $\gP$ is zero and
so $\gF$ is in $\she_2$.

Now we look at the situation when $a \geq 3$. We know then that $\gF$ is 
reflexive. Let
\[ \gF \vpil \gE_0 \vpil \cdots \vpil \gE_m \]
be a locally free resolution of $\gF$. Dualizing we get
\[ \gF^\vee \pil \gE_0{}^\vee \pil \cdots \pil \gE_m{}^\vee \]
with cohomology $\gext^i(\gF, \ompnm)$ at $\gE_i{}^\vee$.
Now we claim that $\gext^i(\gF^\vee, \ompnm)$ is zero for $0 < i < a-1$.
More generally, if $\gK_r$ is the kernel of 
$\gE_r{}^\vee \pil \gE_{r+1}{}^\vee$, then $\gext^i(\gK_r, \ompnm)$ is zero
for $0 < i < r+a-1$. This follows easily by breaking the complex 
$\gE_\cdot$ into short exact sequences and using the fact noted above.
If we take a locally free resolution $\gD_\cdot$ of $\gF^\vee$ and dualize 
this we get a complex
\[ \gF^{\vee \vee} \pil \gD_{a-1}{}^\vee \pil \cdots \pil \gD_0{}^\vee 
\pil \gD_{-1}{}^\vee \pil \cdots \]
which is exact at all $\gD_i{}^\vee$ when $i > 0$. Hence $\gF$ which is 
equal to $\gF^{\vee \vee}$ is in $\she_a$.
\end{proof}

To obtain even nicer classes of coherent sheaves, we can 
consider the subcategory $\she_a^\dagger$ of $\she_a$ consisting of $\gF$ such 
that there exist $r > 0$ with $\gext^i(\gF, \ompnm)$
zero except possibly when $i$ is $0$ or $r$.

\section{Simplicial complexes via the BGG-correspondence}

Given a simplicial complex $\dl$ over $[n]$, recall the module
$C_\dl$ over $E(V)$. Using the BGG-correspondence, we may form the
complex $\tL(C_\dl)$ of coherent sheaves on $\pnm$. Of interest then is the
cohomology sheaves of $\tL(C_\dl)$. The following basic result was 
established in \cite{FV}.

\begin{theorem} a) $\tL(C_\dl)$ has at most one non-zero cohomology 
group $\gF$ if and only if $\dl$ is CLeray. In this case
$\gF = H^{-c} \tL(C_\dl).$

b) The cohomology sheaf $\gF$ is locally Cohen-Macaulay sheaf if and only
if $\dl$ is both CLeray and Cohen-Macaulay.
\end{theorem}

\note We write $L(\dl)$ for $L(C_\dl)$ and $S^a(\dl)$ and $\gS^a(\dl)$ for
the $a$'th cohomology groups of $L(\dl)$ and $\tL(\dl)$ respectively.
In case $\dl$ is CLeray we simply write $\gS(\dl)$ for $\gS^c(\dl)$. If
$\dl$ is both CLeray and Cohen-Macaulay we call $\dl$ bi-Cohen-Macaulay.
\notefin

We shall now establish a more general result, namely determine for any
$\dl$ which cohomology sheaves of $\tL(\dl)$ are non-zero.
We first establish some other results.

\begin{lemma} The associated primes of each cohomology group $H^{-p}L(\dl)$
are of the form $(x_i)_{i \in R}$ where $R$ is a subset of $[n]$.
\end{lemma}

\begin{proof} $H^{-p}L(\dl)$ is graded by $\hele^{n+1}$. By \cite[Exc. 35]{Ei}
an associated prime is also graded by $\hele^{n+1}$. But the only such prime
ideals are of the above form.
\end{proof}

\subsection{Local ranks}
The following Proposition and Corollary is due to T. Ekedahl.
Let $v = \sum \lambda_i e_i$ be a linear form in $E(V)$. It can also be
considered as a point in $\pnm$. We call the set 
$R = \{ i \, | \, \lambda_i \neq 0 \}$ the {\it support} of $v$.

\begin{proposition}
Let $v = \sum \lambda_i e_i$ be a point in $\pnm$ with support $R$.
Then 
\[ H^{-p}(\tL(\dl) \te_{\opnm} k(v)) = \oplus_{S \sus R^c} 
\tvH_{p-1-s}(\dl^{S,T}_R). \]
\end{proposition}

\begin{proof}
First note the following general fact. If $v = \sum \lambda_i e_i$ where
all $\lambda_i$ are non-zero, let $u = \sum_{i=1}^n e_i$. Then 
\[ H^{-p} (C_\dl, v) = H^{-p}(C_\dl,u) = \tvH_{p-1}(\dl).\]
This is because the map $C_\dl \pil C_\dl$ defined
by $x_i \mapsto \lambda_i x_i$ gives an isomorphism between
the complexes $(C_\dl,v)$ and $(C_\dl,u)$.

Now by Paragraph \ref{2Tor} 
\[ H^{-p} (\tL(\dl) \te_{\opnm} k(v)) = H^{-p}(C_\dl,v). \]
Let $V_R$ be the subspace of $V$ spanned by the $e_i$ where $i$ is in $R$.
Restricting $C_\dl$ via $E(V_R) \pil E(V)$ we get
\[ \res \, C_\dl = \oplus_{S \sus R^c} C_{\dl^{S,T}_R}(s). \]
Therefore
\[ H^{-p} (C_\dl, v) = \oplus_{S \sus R^c} H^{-p+s}(C_{\dl^{S,T}_R},v) 
= \oplus_{S \sus R^c} \tvH_{p-s-1}(\dl^{S,T}_R). \]
\end{proof}

\begin{corollary} \label{3rang}
a) The rank of $H^{-p} \tL(\dl)$ is the dimension of $\tvH_{p-1}(\dl)$.

b) Let $\dl$ be CLeray and $v$ a point in $\pnm$ with support $R$.
Then 
\[ \Tor_p^{\opnm}(\gS(\dl), k(v)) = \oplus_{S \sus R^c} 
\tvH_{p+c-1-s} (\dl^{S,T}_R). \]
\end{corollary}

\begin{proof}
Immediate.
\end{proof}

\eks Let the dimension of $\dl$ be $1$, so $\dl$ is a graph and suppose
$c=1$.
Then $\dl$ is CLeray
iff $\dl$ is a forest, $\dl$ is Cohen-Macaulay iff $\dl$ is connected
and $\dl$ is bi-Cohen-Macaulay iff $\dl$ is a tree.

We can compute the rank of $\gS(\dl)$ at a point $v$ with support $R$
as follows. (We use lower case letters to denote the dimension of
a cohomology group.)
\[ \rank \gS(\dl)_v = \tvh_0 (\dl_R) + \sum_{x \in R^c} 
\tvh_{-1} (\dl^{\{x \}}_R). \]
Note that $\tvh_{-1}(\dl^{\{x \}}_R)$ is non-zero iff
$\dl^{\{x \}}_R$ is $\{ \empt \}$ and this holds iff
$R \cap \lkd \{ x \}$ is $\empt$. Thus
\[ \rank \gS(\dl)_v = \tvh_0(\dl_R) + 
\# \{ x \, | \, R \cap \lkd \{x\} = \empt \}.\]
\eksfin

\subsection{Cohomology modules}

Now we shall investigate the cohomology modules of $L(\dl)$. Note that the
complex $L(\dl)$ is simply the {\it cellular complex}, see \cite{BS}, obtained
by associating to vertex $i$ in $\dl$ the monomial variable $x_i$.
For $\bf{a}$ in ${\bf N}^n$ 
let $\dl_{\leq \bf{a}}$ 
be the subcomplex of $\dl$ on the vertices $i$ for which
$a_i > 0$, and call the set of such $i$ the support of $\bf{a}$.
The following observation
is part of the proof of \cite[Prop.1.2]{BS} and gives the multigraded
Hilbert functions of the cohomology modules.

\begin{proposition}
a) $H^{-p} L(\dl)_{\bf a}$ is isomorphic to $\tvH_{p-1}(\dl_{\leq \bf{a}})$.

b) Let $m$ be a monomial and ${\bf b} = {\bf a} + \deg_{{\bf Z}^n} m$.
There is then a commutative diagram
\[ \begin{CD}
\tvH_{p-1}(\dl_{\leq \bf a}) @>>> \tvH_{p-1}(\dl_{\leq \bf b}) \\
@VVV  @VVV \\
H^{-p}L(\dl)_{\bf a} @>{\cdot m}>> H^{-p}L(\dl)_{\bf b}
\end{CD} \]
\end{proposition}

\begin{proof}
The graded part $L(\dl)_{\bf a}$ identifies naturally with the complex 
$(C_{\dl_{\leq \bfa}},d)$ as follows. Let $M$ be the monomial corresponding
to $\bfa$. Now 
\[L^{-p}(\dl) = S(-p) \te_k (C_\dl)_{-p} = \oplus_I S u_I \]
where $u_I = 1 \te_k x_I$. There is a natural map
\[ (C_{\dl_{\leq \bfa}})_p \pil L^{-p}(\dl)_{\bfa}\]
given by
\[ x_I \mapsto M/x_I \cdot u_I. \]
We verify easily that this is an isomorphism and commutes with the 
differentials in $C_{\dl_{\leq \bfa}}$ and $L(\dl)_{\bfa}$. 
This gives a) and b). 
\end{proof}

\begin{corollary} $H^{-p} L(\dl)$ is zero iff $\tvH_{p-1}(\dl_R)$ is 
zero for all $R \sus [n]$.
\end{corollary}

Now we look closer at the cohomology modules of $L(\dl)$. For a subset
$R$ of $[n]$ let $P(R)$ be the homogeneous prime ideal in the polynomial ring
$S$ defining the linear subspace ${\bf P}^{n-1}$ spanned by the $e_i$ where
$i$ is in $R$. That is $P(R) = (x_i)_{i \not \in R}$.

\begin{theorem} \label{4Thm-Ass}
a) The associated primes of $H^{-p}L(\dl)$ are the prime ideals $P(R)$
where $R$ is maximal such that $\tvH_{p-1}(\dl_R)$ is non-zero.

b) The depth of $H^{-p} L(\dl)$ is $\geq p+1$ for $p \geq 1$.

c) For $R$ as in a) there is an isomorphism 
\[ (S/P(R))(-r)^\rho_{(P(R))} \iso H^{-p} L(\dl)_{(P(R))} \]
where $\rho$ is the dimension of $\tvH_{p-1}(\dl_R)$.

\end{theorem}

\begin{proof}
a) Let $R$ be maximal  such that $\tvH_{p-1}(\dl_R)$ is non-zero. Let $m$ in 
$H^{-p} L(\dl)$ have support $R$. Then $x_im = 0$ for $i$ not in $R$,
so clearly $P(R)$ is the annihilator of $m$ and is an associated prime.

Conversely, let $P(R)$ is an associated prime, the annihilator of
a multihomogeneous element $m$. Then $x_i m$ is zero for $i$ not in $R$,
and all $\Pi_{i \in R} x_i^{r_i} m$ are non-zero. 
Hence $H^{-p} L(\dl)_{\bfa}$ is nonzero
if the support of $\bfa$ is $R$ and zero if the support of $\bfa$ strictly 
contains $R$ and so $R$ is maximal such that $\tvH_{p-1}(\dl_R)$ is nonzero.

b) Let 
\[ H^{-p} L(\dl) \vpil F_0 \vpil F_1 \vpil \cdots \vpil F_r \]
be a minimal multigraded resolution. Say that $\bfa$ in ${\bf N}^n$ is of
{\it characteristic} type if every $a_i$ is $0$ or $1$.
By construction of the complex $L(\dl)$, all the minimal generators of 
$H^{-p}L(\dl)$ are of characteristic type
and of degree $\geq p+1$. It is easily seen that the same holds for
the kernel of $F_0 \pil H^{-p} L(\dl)$, and so for all kernels.
So the generators of $F_l$ are of characteristic type and so have degree
$\leq n$. Hence the length $r$ of the resolution is $\leq n-p-1$. 
By the Auslander-Buchsbaum theorem $H^{-p} L(\dl)$ has depth $\geq p+1$.

c) Let $\delta(R)$ in ${\bf N }^n$ have $1$ in position $i$ if $i$ is not
in $R$ and $0$ otherwise. Since $H^{-p} L(\dl)_{\delta(R)}$ is 
isomorphic to $H^{p-1}(\dl_R)$ there is a sequence
\[ 0 \pil (S/P(R))(-r)^\rho \pil H^{-p} L(\dl) \pil Q \pil 0 \]
where $\rho$ is the $k$-dimension of $\tvH_{p-1} (\dl_R)$ and $Q$ is the
cokernel. Now localizing we find that $Q_{(P(R))}$ is zero since
$\Pi_{i \in R} x_i \cdot Q$ is zero.
\end{proof}

\begin{corollary}
a) The sheaf $H^{-p} \tL (\dl)$ has no embedded linear subspaces and
is supported on a reduced union of coordinate linear subspaces of 
${\bf P}^{n-1}$. 

b) For $p \geq 1$, $H^{-p} L(\dl)$ is the graded global
sections of the sheaf $H^{-p} \tL(\dl)$.
\end{corollary}

\begin{proof}
The first statement follow from a) and c) of Theorem \ref{4Thm-Ass}. 
To prove the last statement, let
$M = H^{-p} L(\dl)$. By b) of Theorem \ref{4Thm-Ass} 
the local cohomology groups 
$H_m^0(M)$ and $H_m^1(M)$ are zero, \cite[Thm. A4.3]{Ei}. That
$\oplus_{k\in \hele} \Gamma({\bf P}^{n-1}, \tM(k))$ 
is equal to $M$ now follows by 
\cite[Thm. A4.1]{Ei}.
\end{proof}

\section{Hierarchies of CLeray complexes}

We consider the various hierarchies of coherent sheaves $\she_a$
and subclasses of them and investigate simplicial complexes $\dl$ such that
$\gS(\dl)$ is in such a class.

\subsection{$a+1$-CLeray simplicial complexes}

Let $\cler_a$ be the class of CLeray simplicial complexes $\dl$ such that
$\gS(\dl)$ is in $\she_a$. 
Note that if $\dl$ is the $c-1$-skeleton of the $n-1$-simplex then
$\gS(\dl)$ is the vector bundle $\Omega^c_{{\bf P}^{n-1}}$ and so $\dl$ is
in $\cler_a$ for all $a$. We shall however make the convention 
that the $c-1$-skeleton of the $n-1$-simplex is in $\cler_a$ iff $c \geq a-1$.
Thus for instance $\cler_3$ contains the $1$-skeleton but not the 
$0$-skeleton, $\{ \emptyset \}$, or $\emptyset$.
The simplicial complexes in $\cler_a$ are now called {\it $a+1$-CLeray}. 
The following gives a criterion for $\dl$ to be $a+1$-CLeray.

\begin{theorem} $\dl$ is $a+1$-CLeray iff $c \geq a-1$ and
\[ \tvH_{c-a}(\lk_{\dl} S) = 0, \quad \text{for } s \geq a.\]
In particular, if $\dl$ is not the $a-2$-skeleton of the $n-1$-simplex
then $c \geq a$. \label{4krit-cler}
\end{theorem}

\begin{proof}
 a) That $\dl$ is in $\cler_a$ means that 
$\gext^i(\gS(\dl), \ompnm)$
has codimension $\geq i+a$ for all $i \geq 1$. By (\ref{2extdual})
this sheaf is the cohomology sheaf
\[ H^{c+i-n} \tL((C_\dl)^\vee). \]
Now by (\ref{1check}) there is an exact sequence
\[ 0 \pil C_{\dl^*} \pil E(W) \pil (C_\dl)^\vee \pil 0 \] 
and so we get
$\gext^i(\gS(\dl), \ompnm)$ equal to $H^{c+i-n+1} \tL(\dl^*)$. That this is a
sheaf with support in codimension $\geq i+a$ means that 
\begin{equation}\tvH_{n-i-2-c}((\dl^*)_{-S}) = 0 \label{4.1} 
\end{equation}
when $s < i+a$. By Lemma \ref{1alexdual} $(\dl^*)_{-S}$ is $(\lkd S)^*$. 
And so by Lemma \ref{1alexdualhom}, (\ref{4.1})
is the same as
\begin{equation} \tvH_{c+i-s-1}(\lkd S) = 0 \label{4.2}
\end{equation}
for $s < i+a$. Now put $s =  i+a-1-p$ where $p \geq 0$. Then this becomes
\begin{equation} 
\tvH_{c-a+p}(\lkd S) = 0
\label{4.3}
\end{equation}
for $ s \geq a-p$.
By Corollary \ref{1redcor} c) 
if $c \geq a$ then 
$\tvH_{c-a} (\lkd S)$ zero 
for $s \geq a$ implies (\ref{4.3}), so this latter is 
the condition that $\dl$ is in $a+1$-CLeray when $c \geq a$. When 
$c < a$, then again by Corollary \ref{1redcor} c), (\ref{4.3}) is implied
by $\tvH_{-1}(\lk_{\dl} S)$ zero for $s \geq c+1$ and this means that
$\dl$ has no faces of dimension $c$. Hence it is the $c-1$-skeleton of 
the $n-1$-simplex and by our convention we have $c = a-1$.
\end{proof}

\begin{corollary} \label{4-lkcor}
a) $\dl$ is $a+1$-CLeray iff $\tvH_p(\dl^S_R)$ is zero for $p+s \geq c$
and $p \geq c-a$.

b) If $\dl$ is $a+1$-CLeray, then $\dl_R$ is $a+1$-CLeray.
\end{corollary}

\begin{proof}
a) Suppose $\tvH_p(\dl^S_R)$ is nonzero where $p+s \geq c$ and $ p \geq c-a$
where $c \geq a$. By Corollary \ref{1redcor} e) $H_p(\lk_\dl S^\prime)$
is nonzero for some $S^\prime$ containing $S$ and by Corollary \ref{1redcor} c)
we get  $\tvH_{c-a}(\lk_\dl S^{\prime \prime})$ is nonzero where
$s^{\prime \prime} - s^\prime \geq p -(c-a)$. But then
\[ s^{\prime \prime} \geq p+s^\prime - (c-a) \geq c-(c-a) = a \]
and this is against assumption.

b) This is clear since the frame dimension of $\dl_R$ is greater
or equal to that of $\dl$.
\end{proof}

As stated after Definition \ref{1cler}, when $\dl$ is $1$-CLeray
with $c=0$, then $\dl$ is a simplex on its vertices. In view of 
Theorem \ref{4krit-cler} it is of interest to investigate the 
$a+1$-CLeray $\dl$'s with minimal interesting frame dimension which is $a-1$
(when the frame dimension is $a-2$ it is the $a-2$-skeleton of the 
$n-1$-simplex).

\begin{proposition} \label{4prop-c+1CLeray}
$\dl$ is $c+1$-Cleray iff any two
distinct facets intersect in a subset of cardinality less than $c$.
(If $c=0$ this means that there is only one facet.)
\end{proposition}

\begin{proof}
Suppose $\dl$ is $c+1$-CLeray. Let $F$ and $G$ be two distinct facets such 
that the cardinality of $S = F \cap G$ is maximal. Now we claim that $\lkd S$
is disconnected. Suppose not and let $F = F_1 \cup S$ and $G = G_1 \cup S$
where $F_1 \cap G_1$ is empty. Then there would be a path from some vertex 
$f_1$ in $F_1$ to some vertex $g_1$ in $G_1$ in $\lkd S$. Say the path starts
with $\{f_1,x\}$ where $x$ is not a vertex in $F_1$. 
Then $H = S \cup \{f_1, x\}$ is
a face of $\dl$, $H$ is not contained in $F$ and $F \cap H$ has cardinality
larger than $S$ which goes against our assumptions. Hence $\tvH_0 (\lkd S)$ is
non-zero. By Theorem \ref{4krit-cler} a) we get that $s < c$.
 
In the other direction, given that distinct 
facets always intersect in cardinality
$< c$, we see that $\lkd S$ when $s \geq c$ is always a simplex. Hence
$\tvH_0 (\lkd S)$ is zero and so $\dl$ is $c+1$-Cleray.
\end{proof}

By Corollary \ref{1redcor}, if $p+s \leq c-2$ then $\tvH_p(\lk_{\dl} S)$
is zero and Corollary \ref{4-lkcor} gives conditions on homology 
when $p+s \geq c$. It is therefore
of interest to investigate what happens when $p+s = c-1$. The following
is motivating for the classes studied in Subsection \ref{4cler-ring}.

\begin{proposition} ($p+s = c-1$) 

a) If $\tvH_{c-s-1} (\lk_{\dl} S)$ is zero and $S^\prime$ contains $S$ then
$\tvH_{c-s^\prime -1}(\lk_{\dl} S^\prime)$ is zero.

b) Suppose $\dl$ is $a+1$-CLeray. Then $\tvH_{c-a}(\lk_{\dl} S)$ is nonzero
for any face $S$ with $s = a-1$. \label{4-cutprop}
\end{proposition}

\begin{proof}
 a) follows form Corollary \ref{1redcor} a) by contraposition.

b) Let $S^\prime$ be a facet containing $S$ where $s$ is $a-1$ 
Then $\tvH_{-1}(\dl^{S^\prime,\empt}_{-S^\prime})$ is non-zero. Reduce 
$S^\prime$ to $S$ and get $\tvH_p (\dl^{S,T}_R)$ non-zero and so
$p+s \geq c-1$. 
Now reduce $T$ to $\empt$ and get $\tvH_p(\lkd \tS)$ non-zero
where $\tS \sups S$.
Since $\dl$ is in $\cler_a$, $\tvH_p (\lkd \tS)$ is zero when $p + \ts \geq c$
and $p \geq c-a$. Hence either $p+\ts \leq c-1$ or $p \leq c-a-1$. The
latter is impossible since $p+a-1$ which is $p+s$ is $\geq c-1$
and the former gives $\ts = s = a-1$ and so
$\tvH_p( \lkd S)$ is non-zero.
\end{proof}

We now give a somewhat more conceptual description of what it means
for a complex to be $a+1$-CLeray. First a lemma.

\begin{lemma} \label{4-2CLeray}
If $\dl$ is $2$-CLeray, the frame dimension of $\lk_{\dl} \{x\}$ is one
less than that of $\dl$. 
\end{lemma}

\begin{proof} Since $\dl$ is 2-CLeray, $c \geq 0$. 
If $c= 0$ this holds, so assume $c \geq 1$.
We have $\tvH_p (\lk_\dl( \{ x \} \cup S))$ zero for $p+s+1 \geq c$ and
$p \geq c-1$, which reduces simply to the condition $p \geq c-1$. 
Thus by Lemma \ref{1lem-clk} the frame dimension of $\lk_{\dl} \{x \}$
is $\leq c-2$ and so must be equal to $c-2$.
\end{proof}

\begin{theorem}
$\dl$ is $a+1$-CLeray iff every link $\lk_{\dl} S$ with $s=a$ is 
CLeray with frame dimension $c-a-1$. \label{4-lkthm}
\end{theorem}

\begin{proof} We may assume $c \geq a \geq 1$. If $\dl$ is $a+1$-CLeray
then 
\[\tvH_{c-a} \lk_\dl (S \cup T)= 0, 
\quad \text{ when } s=a \text{ and } t \geq 0. 
\label{4SUT}\]
Since by Lemma \ref{4-2CLeray} the frame dimension of $\lk_\dl S$ is
$c-a-1$, we get that $\lk_\dl S$ is CLeray. Conversely, if each
$\lk_\dl S$ is CLeray of frame dimension $c-a-1$ then (\ref{4SUT})
holds and so $\dl$ is $a+1$-CLeray.
\end{proof}

\subsection{ The classes $\cler_a^\dagger$}

Let $\clerd_a$ be the class of CLeray simpicial complexes $\dl$ such that
$\gS(\dl)$ is in $\she_a^\dagger$. The following gives a criterion for 
$\dl$ to be in $\clerd_a$.

\begin{theorem}  $\dl$ is in $\clerd_a$ iff $\tvH_p(\lk_\dl S)$ is zero 
for all $p$ and $s$ in the range
\[ c \leq p+s \leq d-2 \]
and also when $p,s$ are in the range 
\[ p+s = d-1,\,\, p \geq c-a. \]
In this case every facet has dimension $c-1$ or $d-1$. \label{4-dthm}
\end{theorem}

\begin{proof}
For $\gS(\dl)$ to be in $\cler_a^\dagger$, there is some $r$ such that 
$\gext^i(\gS(\dl), \ompnm)$ is zero except possibly when $i$ is $0$ or $r$, 
and in the latter case it has codimension $\geq r+a$. By the argument 
of Theorem \ref{4krit-cler} this means that 
\begin{eqnarray} &H_{c+i-s-1} (\lkd S) = 0 \quad & \text{for all $s$ when } 
i \neq 0, r \label{4.4} \\
&H_{c+r-s-1} (\lkd S) = 0 \quad & \text{for } s < a+r. \notag
\end{eqnarray}

Letting $p = c+i-s-1$ this gives $H_p(\lkd S)$ zero when 
$p+s \geq c$ and $p+s \neq c+r-1$ or when $p+s = c+r-1$ and $p \geq c-a$.
But let $S$ be a facet. Then $\lkd S$ is $\{\empt \}$ and so 
$\tvH_{-1}(\lkd S)$ is non-zero. Hence if $s \geq c+1$, letting $s = c+i$,
we see from (\ref{4.4}) that we must have $i=r$ and so $s=c+r$. 
So all facets $S$ with $s \geq c+1$ must have $s = c+r$ 
which is then equal to $d$. 
\end{proof}

\begin{corollary} \label{4-dcor} 
a) $\dl$ is in $\clerd_a$ iff $\tvH_p (\dl^S_{-R})$ is zero when 
$p,s$ (and $r$)
are in the range
\[ c \leq p+s \leq d-2 \]
and also when $p$ and $s$ are in the range
\[p+s = d-1, \,\, p \geq c-a. \]

b) If $\dl$ is in $\clerd_a$ then $\dl_{-R}$ is in $\clerd_a$.
\end{corollary}

\begin{proof} The proofs are analogous to those of Corollary \ref{4-lkcor}.
\end{proof}

\begin{theorem} \label{4-dthm2}
$\dl$ is in $\clerd_a$ iff every link $\lk_\dl S$ with $s=a$
is in $\clerd_0$ with dimension $d-a-1$ and frame dimension $c-a-1$ 
or is the $c-a-1$ skeleton of the simplex on $[n]\backslash S$.
\end{theorem}

\begin{proof}
We can assume $c \geq a \geq 1$. If now $\dl$ is in $\clerd_a$, the
facets have dimension $c-1$ or $d-1$. In view of Lemma \ref{1lem-clk}
the statement
follows immediately from Theorem \ref{4-lkthm}.
\end{proof}

\rem The if direction is not true unless we assume the links to be
of dimension $d-a-1$ or $c-a-1$. If one simply assume the links are of 
frame dimension $c-a-1$ and are in $\clerd_0$, then a counterexample is 
given by starting with the disjoint union of a $2$-simplex and a $3$-simplex
and then adding all line segments between pairs of vertices of them.
\remfin

\subsection{ The classes $\clerri_a$} \label{4cler-ring}

By Proposition \ref{4-cutprop} 
the nicest behaviour one can expect for the homology groups
$\tvH_p(\lk_\dl S)$ when $p+s = c-1$ and $\dl$ is in $\cler_a$ is 
\begin{equation}
 \tvH_{c-a-1} (\lk_\dl S) = 0, \quad \text{when } s=a. \label{4cut}
\end{equation}
We now let $\clerri_a$ be the complexes $\dl$ such that $\dl$ is in 
$\clerd_a$ and fulfills the condition (\ref{4cut}).

(Note that when $\dl$ is in $\cler_{a+1}$ then by Proposition \ref{4-cutprop}, 
$\tvH_{c-a-1}(\lk_\dl S)$ is nonzero for $s=a$, so if $\dl$ is in $\cler_a$
and fulfills (\ref{4cut}) then it is not in $\cler_{a+1}$.)

Two special cases are interesting to take note of.

\eks If $\dl$ in $\clerri_a$ has $c=a$, the lowest interesting value
for $c$, then the condition (\ref{4cut}) says that $H_{-1} (\lk_\dl S)$
is zero for $s=c=a$. Hence by Theorem \ref{4-dthm} all facets of $\dl$ have
dimension $d-1$, and by Proposition \ref{4prop-c+1CLeray} 
any two facets intersect in a face
of dimension $\leq c-2$. Hence we get precisely the Steiner systems
$S(c,d,n)$. In particular the $f$-vector only depends on $c,d$, and $n$. 
\eksfin

\eks When $\dl$ is in $\clerri_0$ the condition (\ref{4cut}) says 
$\tvH_{c-1} (\dl)$ is zero and so $\gS(\dl)$ is a torsion sheaf.
Then $\gext^i(\gS(\dl), \omega_{{\bf P}^{n-1}})$ is nonzero only for $i = d-c$
and so $\gS(\dl)$ is a locally Cohen-Macaulay sheaf.
By \cite{FV} such $\dl$ are bi-Cohen-Macaulay, i.e. both $\dl$ and
$\dl^*$ are Cohen-Macaulay and we showed that the $f$-vector of such $\dl$
again only depends on $n,d$, and $c$.
In fact we showed that the $f$-polynomial $f_\dl(t) = \sum_i f_{i-1} t^i$
is 
\[ f_\dl(t) = (1+t)^{d-c}(1+(n-d+c)t + \cdots + 
\binom{n-d+c}{c} t^c). \]
\eksfin

By the following theorem and its corollary this generalizes.

\begin{theorem} \label{4-rithm}
$\dl$ is in $\clerri_a$ iff each $\lk_\dl S$ with $s=a$ 
is bi-Cohen-Macaulay with invariants $n-a$, $d-a$, and $c-a$.
\end{theorem}

\begin{proof}
Suppose $\dl$ is in $\clerri_a$. The condition (\ref{4cut}) together
with Proposition \ref{4-cutprop} a) gives $H_{-1}(\lk_\dl S)$ zero when $s=c$. 
Hence by Proposition
\ref{4prop-c+1CLeray} all facets of $\dl$ have dimension $d-1$. 
This shows that $\lk_\dl S$ has
invariants $n-a, d-a$, and $c-a$. It is now immediate that $\lk_\dl S$
fulfills the condition to be in $\clerri_0$ and this means that 
$\lk_\dl S$ is bi-Cohen-Macaulay. In the converse direction it is immediate
to see that $\dl$ fulfills the condition to be in $\clerri_a$.
\end{proof}

\begin{corollary} Suppose $\dl$ in $\clerri_a$
has invariants $n,d$, and $c \geq a$. 
Then $f_\dl(t)$ is a 
polynomial $f_a(n,d,c;t)$ depending only on $a,n,d$, and $c$, and is given
inductively as follows.

1) $f_a^\prime(n,d,c;t)/n = f_{a-1}(n-1,d-1,c-1; t)$

2) $f_0(n,d,c;t) = (1+t)^{d-c}(1+(n-d+c)t + \cdots + 
\binom{n-d+c}{c} t^c).$
\end{corollary}

\begin{proof} We may assume $a \geq 1$ and assume
first that $c \geq 2$. By induction the $f$-polynomial of $\lkd \{x\}$
is 
\[ g(t) = f_{a-1}(n-1,d-1,c-1;t) \]
and is independent of $x$. 
Now if $F$ is any $i$-dimensional face of $\dl$, then
$F\backslash \{x \}$ is an $i-1$-dimensional face of $\lkd \{x\}$ 
for each $x$ in $F$.
Therefore $n g_{i-1}$ counts all pairs $(F,x)$ where $x$ is in $F$. 
But this is also equal to to $f_i (i+1)$ and so 
\begin{eqnarray*}
ng_{i-1} & = & f_i (i+1) \\
n g(t) & = & f_\dl^\prime (t).
\end{eqnarray*}

In case $c=1$ then $a=1$ and $\dl$ is a disjoint union of $n/d$ simplexes
of dimension $d-1$ and so 
\begin{equation*} f_\dl(t)  =  (n/d) (1+t)^d + 1 - n/d
\end{equation*}
and so 1) also holds in this case since 
\begin{equation*}
f_0(n-1, d-1, 0;t) =  (1+t)^{d-1},
\end{equation*}
\end{proof}

\rem The $\dl$ in $\clerri_a$ are
$a$-$(n,d,\lambda)$ block designs (see \cite[Chap.14]{Br}), where $\lambda$
is the number of facets of $\lkd S$ for $S$ any face of cadinality
$a$. One may see that $\lambda$ is $\binom{n-d+c-a}{c-a}$.
\remfin

The following generalizes the fact that when $\dl$ is in $\clerri_0$ 
then $\gS(\dl)$ is a locally Cohen-Macaulay sheaf of codimenison $d-c$ and so 
$\gext^{d-c}(\gS(\dl), \ompnm)$ is a locally Cohen-Macaulay sheaf.

\begin{proposition}
If $\dl$ is in $\clerri_a$, then the only 
non-vanishing higher $\gext$-sheaf, $\gext^{d-c}(\gS(\dl), \ompnm)$, is
a locally Cohen-Macaulay sheaf.
\end{proposition}

\begin{proof}
For short write $\gE$ for $\gext^{d-c}(\gS(\dl), \ompnm)$.
The sequence
\[ 0 \pil C_{\dl^*} \pil E(W) \pil (C_\dl)^\vee \pil 0 \]
gives an exact sequence
\[ 0 \pil \tL(\dl^*) \pil \tL(E(W)) \pil \tL((C_\dl)^\vee) \pil 0.\]
Hence by (\ref{2dualkx}) and (\ref{2extdual}), 
$\tL(\dl^*)$ only has possibly non-zero cohomology
sheaves
$\gE$ in degree $c+ (d-c)-n+1$ which is $-c^*$,
and $H^{-d^*}\tL(\dl^*)$ in degree $-d^*$. 
The condition (\ref{4cut}) for $\dl$, when translated to $\dl^*$, see 
(\ref{5cut}) and the paranthetical remark after, gives that all generators
of $H^{-d^*}L(\dl^*)$ have degree $\geq n-a+1$. Since they are of 
characteristic type (see the proof of Theorem \ref{4Thm-Ass}) a free
resolution has length $\leq a-1$.


Hence $\gE$ has a locally free resolution of length
$\leq d^*-c^*+a$ which is equal to $d-c+a$.
By the Auslander-Buchsbaum theorem \cite[Thm. 19.9]{Ei},
$\gE$ then has local depth $\geq (n-1) - (d-c+a)$.
But $\dl$ being in $\cler_a$, this sheaf has local dimension 
$\leq (n-1) - (d-c+a)$.
Since local dimension is greater or equal to local depth we must have
equalities everywhere and so $\gE$ is 
locally Cohen-Macaulay.
\end{proof}

\section{Hierarchies of Cohen-Macaulay simplicial complexes}

The Alexander duals of CLeray simplicial complexes are the Cohen-Macaulay
simplicial complexes. Therefore by taking Alexander duals of the various
hierarchies of CLeray complexes, we get hierarchies of Cohen-Macaulay
simplicial complexes. This section will contain the Alexander dual versions
of most of the statements in Section 4. But we shall also give a description
of what it means for a simplicial complex $\dl$ to be Gorenstein* in 
terms of the associated coherent sheaf $\gS(\dl^*)$. We prove that
$\dl$ is Gorenstein* iff $\gS(\dl^*)$ is an ideal sheaf, i.e. a subsheaf
of $\go_{{\bf P}^{n-1}}$. In fact it turns out to be the associated sheaf
of the ideal defining the Stanley-Reisner ring of $\dl^*$.

\subsection{$a+1$-Cohen-Macaulay simplicial complexes}
Let $\cm_a$ be the class of Cohen-Macaulay simplicial complexes which are
Alexander duals of the simplicial complexes in the class $\cler_a$.
We shall show that this is exactly the class of {\it $a+1$-Cohen-Macaulay} 
simplicial complexes as defined by Baclawski \cite{Ba}.

\begin{theorem} $\dl$ is in $\cm_a$ iff $d \leq n-a$ and 
$\tvH_p(\dl_{-R})$ is zero when $p+r = d+a-2$ and $p \leq d-2$.
In particular, if $\dl$ is not the $n-a-1$-skeleton of the $n-1$-simplex,
then $d \leq n-a-1$. \label{5thm-hoved}
\end{theorem}

\begin{proof}
 By Theorem \ref{4krit-cler} $\dl^*$ is in $\cler_a$ means that $d \leq
n-a$ and 
\begin{equation} \tvH_{c^*-a}((\dl^*)^S_{-S}) = 0 \quad \text{for } s \geq a.
\label{5.1}
\end{equation}
Now $(\dl^*)^S_{-S}$ is equal to $(\dl_{-S})^*$ and so by Lemma 
\ref{1alexdualhom},
(\ref{5.1}) is equivalent to 
\[ \tvH_{(n-s)+a-c^*-3}(\dl_{-S}) = 0 \quad \text{for } s \geq a. \]
Since $n-1-c^* = d$ we get the statement.
\end{proof}

\begin{corollary} \label{5cor-range}
a) $\dl$ is in $\cm_a$ iff $\tvH_p(\dl^S_{-R})$
is zero when $p+s \leq d-2$ and $p+r \leq d+a-2$.

b) If $\dl$ is in $\cm_a$, then $\lkd S$ is also in $\cm_a$.
\end{corollary}

\begin{proof}
This is just the Alexander dual versions of Corollary \ref{4-lkcor},
using Lemma \ref{1alexdualhom}.
\end{proof}

\begin{corollary}
Suppose $\dl$ has dimension one, i.e. $\dl$ is a graph.
Then $\dl$ is in $\cm_a$ if and only if $\dl$ contains at least $a+2$
vertices and is $a+1$-connected.
\end{corollary}

\begin{proof} $\dl$ is in $\cm_a$ iff $\tvH_{-1} (\dl_{-R})$ is zero for
$r = a+1$ and $\tvH_0 (\dl_{-R})$ is zero for $r =a$. 
This translates precisely to the above.
\end{proof}

The following describes the objects of $\cm_a$
with upper extremal values for $d$. Recall that a subset $F$ of $[n]$ is a 
{\it missing face} of $\dl$ if 
$F$ is not in $\dl$. 
Now if $\dl$ is in $\cm_a$ and not the $n-a-1$-skeleton of the $n-1$-simplex
then $d \leq n-a-1$ or equivalently $a \leq n-d-1$.

\begin{proposition}  $\dl$ is in $\cm_{n-d-1}$
iff the cardinality of $F \cup G$ is $\geq d+2$ for all
distinct missing faces $F$ and $G$.
\end{proposition}

\begin{proof} $\dl$ is in $\cm_{n-d-1}$ iff $\dl^*$ is $c^*+1$-CLeray.
If $F$ and $G$ are two distinct 
missing faces then the complements $F^c$ and $G^c$ are faces of 
$\dl^*$ and so the cardinality of $F^c \cap G^c$ is $\leq c^* - 1$. 
But then the cardinality of $F \cup G$ is $\geq n-c^* +1$ and this is $d+2$. 
\end{proof}

The following shows that $\dl$ is in $\cm_a$ iff $\dl$ is $a+1$-Cohen-Macaulay
($a+1$-CM) as defined by Baclawski, \cite{Ba}.

\begin{theorem}  $\dl$ is in $\cm_a$ iff every restriction 
$\dl_{-R }$ with $r = a$ is Cohen-Macaulay of the same dimension as $\dl$.
\end{theorem}

\begin{proof}
This is just the Alexander dual version of Theorem \ref{4-lkthm}.
\end{proof}

From \cite{Ba} we have the following means of constructing $a+1$-CM
simplicial complexes.

\begin{theorem}[\cite{Ba}] If $\dl$ and $\dl^\prime$ are $a+1$-CM 
simplicial complexes, then the join $\dl * \dl^\prime$ is also $a+1$-CM.
\end{theorem}

In particular, the $\dl^\prime$ consisting of $m \geq a+1$ vertices is 
$a+1$-CM and so the $m$-point suspension of $\dl$
will be $a+1$-CM.

\rem If $F_\cdot$ is a free resolution of the Stanley-Reisner ring of 
$\dl$, it follows from Theorem \ref{5thm-hoved} and Hochsters description
of $\text{Tor}_i^S(k[\dl],k)$, see \cite[II.4.8]{St}, that $\dl$ is
$(a+1)$-CM iff it is Cohen-Macaulay and $F_i = S(-i-d)^{r_i}$ for 
$n-d \geq i > n-d-a$.
\remfin

\subsection{The classes $\cmd_a$}

Let $\cmd_a$ be the class of simplicial complexes which are Alexander
duals of the simplicial complexes in $\clerd_a$.

\begin{theorem} $\dl$ is in $\cmd_a$ iff $\tvH_p(\dl_{-R})$ is zero for
all $r$ when $p$ is in the range 
$c \leq p \leq d-2$ and also when $p = c-1$ and $r < d+a-c$. \label{5-dthm}
\end{theorem}

\begin{proof}
This is the Alexander dual version of Theorem \ref{4-dthm}.
\end{proof}

\begin{corollary}
a) $\dl$ is in $\cmd_a$ iff $\tvH_p(\dl^S_{-R})$ is zero for $p,s$ (and $r$)
in the range 
$ c \leq p+s \leq d-2$
and also when $p+s = c-1$ and $s+r < d+a-c$.

b) If $\dl$ is in $\cmd_a$ then $\lkd S$ is in $\cmd_a$.
\end{corollary}

\begin{proof}
This is the Alexander dual version of Corollary \ref{4-dcor}.
\end{proof}

\begin{theorem}
$\dl$ is in $\cmd_a$ iff each restriction $\dl_{-R}$ with $r = a$
is either in $\cmd_0$ with the same dimension and frame dimension as $\dl$
or is the $d-1$-skeleton of the simplex on $[n]\backslash R$.
\end{theorem}

\begin{proof}
This is the Alexander dual version of Theorem \ref{4-dthm2}.
\end{proof}

\subsection{The classes $\cmri_a$.} \label{5-SekKutt}

We let $\cmri_a$ be the class of simplicial complexes which are Alexander
duals of the simplicial complexes in $\clerri_a$. If $\dl$ is in $\cmri_a$
so $\dl^*$ is in $\clerri_a$ the condition (\ref{4cut}) is equivalent to 
the condition 
\begin{equation} 
\tvH_{d-1}(\dl_{-R}) = 0 \quad \text{when } r=a. \label{5cut}
\end{equation}
(It then follows by Proposition \ref{4-cutprop} that $\tvH_{d-1}(\dl_{-R}) = 0$
for $r \geq a$.) Thus $\dl$ is in $\cmri_a$ iff it is in $\cmd_a$ and fulfills
(\ref{5cut}).

\eks  If $\dl$ is a cyclic polytope of odd dimension, 
then $\dl$ is in $\cmri_1$. To see this, note that
by Alexander duality for Gorenstein* simplicial complexes \cite[p.66]{St}
\begin{equation} 
\tvH_p(\dl_{-R}) = \tvH_{d-2-p}(\dl_R)^* \label{5AD}.
\end{equation}
Now we can apply Theorem \ref{5-dthm}. For a cylic polytope
of odd dimension, $d = 2c$ and we see that if $c \leq p \leq d-2$ then
$0 \leq d-2-p \leq c-2$ and so (\ref{5AD}) is zero.
Also $p = c-1$ gives $d-2-p = c-1$ and $r \leq d+a-c-1$ is the same as 
$r \leq c$. Thus (\ref{5AD}) is zero and so $\dl$ is in $\cmd_1$. That 
(\ref{5cut}) holds is immediate from (\ref{5AD}) and so $\dl$ is in $\cmri_1$.

\begin{theorem}
$\dl$ is in $\cmri_a$ iff every restriction $\dl_{-R}$ where $r=a$,
is bi-Cohen-Macaulay of the same dimension and frame dimension as $\dl$.
\end{theorem}

\begin{proof}
This is the Alexander dual version of Theorem \ref{4-rithm}.
\end{proof}

\begin{corollary}
If $\dl$ is in $\cmri_a$ and $r \leq a$, then $\dl_{-R}$ is in 
$\cmri_{a-r}$.
\end{corollary}

\subsection{Gorenstein* simplicial complexes}

We shall now describe where Gorenstein* complexes fit in our scheme.
Recall that $\dl$ is Gorenstein* if $\tvH_p (\lkd S)$ is $k$ when 
$p+s = d-1$ and $S$ a face of $\dl$, and zero otherwise.
The following theorem is certainly well known (except the last statement). 
It follows for instance from \cite[4.7]{Ba}.

\begin{theorem} $\dl$ is Gorenstein* iff $\dl$ is $2$-CM and 
$\tvH_{d-1}(\dl)=k$, i.e. iff $\gS(\dl^*)$ is a torsion free 
rank one sheaf.  
\end{theorem}

\begin{proof} Suppose $\dl$ is Gorenstein*. 
By Theorem \ref{5thm-hoved} we must show that $\tvH_p (\dl_{-R})$ is zero when
$p+r = d-1$ and $p \leq d-2$. By the Alexander duality theorem
for Gorenstein* simplicial complexes \cite[p.66]{St},
$\tvH_p (\dl_{-R})$ is equal to $\tvH_{d-2-p} (\dl_R)^*$. Hence
we need to show that $\tvH_{p^\prime}(\dl_R)$ is zero when 
$r = p^\prime +1$. But this is true simply because 
the homology of any simplicial complex on $n$ elements always vanishes
in degrees $\geq n-1$.

\medskip

Now suppose $\dl$ is in $\cm_1$ with $\tvH_{d-1} (\dl)$ equal to $k$.
By Corollary \ref{5cor-range} a) then 
$\tvH_p (\lkd S)$ is zero for $p+s \leq d-2$.
We must therefore show in addition that $\tvH_p (\lkd S)$ is $k$ 
when $p+s = d-1$ and $S$ a face of $\dl$. Since $\lkd \{ x \}$ is in $\cm_1$
it is enough by induction to show that $\tvH_{d-2} (\lkd \{ x \})$ is $k$
when $\{ x \}$ is a face of  $\dl$. But there is a sequence
\[ \tvH_{d-1} (\dl_{-\{ x\}}) \pil \tvH_{d-1} (\dl) \pil \tvH_{d-2} 
(\lkd {\{ x \}}) \pil \tvH_{d-2} (\dl_{-\{ x \}}). \]
Since $\dl$ is in $\cm_1$, $\tvH_{d-2} (\dl_{-\{ x \}})$ is zero.
Since $\lkd {\{ x \}}$ is in $\cm_1$, applying Proposition \ref{4-cutprop} b)
we get that $\tvH_{d-2} (\lkd {\{ x \}})$ is non-zero. And so we must
have $\tvH_{d-2} (\lkd {\{ x \}})$ equal to $k$.
\end{proof}

\begin{proposition} \label{5-gorid}
When $\dl$ is Gorenstein* there is a natural identification of 
$S(\dl^*)$ with the ideal defining the Stanley-Reisner ring of the simplicial 
complex $\dl^*$. (And so this ideal defines a subscheme of codimension
$c+1$ in $\pnm$.)
\end{proposition}

\begin{proof}
The paranthetical remark is because the codimension of the subscheme 
defined by the ideal of the Stanley-Reisner ring of $\dl^*$ is $n-d^*$
which is $c+1$.

First we shall suppose $c \geq 1$. Then $\tvH_{d-1} (\dl_{-\{ x \}})$
by Alexander duality is equal to $\tvH_{-1} (\dl_{\{ x \}})$ which is zero.
Hence $\dl$ fulfills condition (\ref{5cut}) and by the last remark in 
Subsection \ref{5-SekKutt}, $H^{-d} \tL(\dl)$ is equal to $\opnm(-n)$.

Now by (\ref{2extdual}) and the sequence
\begin{equation*}
 0 \pil C_\dl \pil E(W) \pil (C_{\dl^*})^\vee \pil 0
\end{equation*}
we get
\[ \ghom(\gS(\dl^*), \ompnm) = H^{-d}\tL(\dl) = \opnm(-n) = \ompnm. \]
Since $\gS(\dl^*)$ is torsion free 
\[ \gS(\dl^*) \inpil \gS(\dl^*)^{\vee \vee} = 
\ghom( \ompnm, \ompnm) = \opnm. \]
Taking graded global sections
\[ S(\dl^*) = \oplus_{m \in \hele}\Gamma(\pnm, \gS(\dl^*)(m)) 
\inpil \oplus_{m \in \hele} \Gamma(\pnm, \opnm(m)) = S\]
we get $S(\dl^*)$ as an ideal in $S$.

If $c= 0$, let $R$ consist of the $x$ such that $\{ x \}$ is a face in
$\dl$. Let $r$ be the cardinality of $R$. 
In this case $H^{-d}\tL(\dl)$ is $\opnm(-r)$. Hence
$\gS(\dl^*)^{\vee \vee}$ is  $\opnm(r-n)$.
By the inclusion 
\[ \opnm(r-n) \mto{\Pi_ {i \not \in R} x_i } \opnm \]
we also get $\gS(\dl^*)$ as an ideal sheaf in $\opnm$ and
taking graded global sections we get $S(\dl^*)$ as an ideal in $S$.

\medskip

Now we want to identify this ideal as the ideal of the Stanley-Reisner ring
associated to $\dl^*$.
By the sequence 
\begin{equation*}
 0 \pil C_{\dl^*} \pil E(W) \pil (C_{\dl})^\vee \pil 0 \label{5bassek}
\end{equation*}
 we get (using (\ref{2dualkx}))
a sequence of complexes
\begin{equation} 
0 \pil L(\dl^*) \pil L(E(W)) \pil \Hom_S(L(\dl), S(-n))[n] \pil 0. 
\label{5Lkx}
\end{equation}

Now the latter complex is
\begin{equation} 
S(d-n)^{f_{d-1}} \vpil S(d-1-n)^{f_{d-2}} \vpil \cdots \vpil S(-n) 
\label{5res}
\end{equation}
where the cohomological degree of the first term is $d-n$. Taking the long
exact cohomology sequence of (\ref{5Lkx}) we get that the only 
cohomology of (\ref{5res}) is in cohomological degree $d-n$ and is
\[ H^{d-n+1} L(\dl^*) = H^{-c^*}L(\dl^*) = S(\dl^*). \]
Hence (\ref{5res}) is a resolution of the ideal $S(\dl^*)$ in $S$.
Since it is a multigraded resolution the generators of $S(\dl^*)$
in $S$ will have multidegrees the multidegrees of the indices  $J$
when writing \[ S(d-n)^{f_{d-1}} = \oplus_J S(d-n) u_J. \] 
But by the BGG-correspondence we recognise the $J$'s as
the multidegrees of the facets of $\dl$. But this means 
exactly that $S(\dl^*)$ is the ideal of the Stanley-Reisner ring
of $\dl^*$.
\end{proof}

\rem By the result of \cite{ER98} the Stanley-Reisner ideal $I_{\dl^*}$ in 
$S$ of a simplicial complex $\dl^*$ has a linear resolution iff $\dl$ is
Cohen-Macaulay. The same is valid for the the analog $J_{\dl^*}$ 
of the Stanley-Reisner ideal in the exterior algebra $E$. But this is the
same as saying that the ideal $J_{\dl^*}$ is a {\it Koszul module} over $E$
(see \cite{Gr}). The algebras $E$ and $S$ are Koszul duals. Via the functors
relating (complexes of) modules over them, $J_{\dl^*}$ then transfers to a 
Koszul module over $S$. The above Proposition \ref{5-gorid} then shows
that this transformed module is an ideal in $S$ exactly when $\dl$ is
Gorenstein*.

Summing up, $\dl$ is Cohen-Macaulay iff $J_{\dl^*}$ in $E$ is a Koszul ideal
(which transfers to a Koszul module), and $\dl$ is Gorenstein* iff $J_{\dl^*}$
in $E$ is a Koszul ideal transforming to a Koszul {\it ideal} in $S$.
\remfin

\rem If $\dl$ is in $\cm_a$ then the link $\lkd S$ is also in $\cm_a$.
Hence if $S$ is of dimension $d-2$, then $\lkd S$ consists of a set of
vertices of, being in $\cm_a$, cardinality $\geq a+1$.

The subclass $\gor_a$ of $\cm_a$ such that this cardinality is always the 
minimum possible, namely $a+1$, might be a reasonable generalization of 
Gorenstein* complexes, since $\gor_1$ would be exactly this class.
If $\dl$ and $\dl^\prime$ are in $\gor_a$ then the join $\dl*\dl^\prime$
is also in $\gor_a$, in particular the $a+1$-point suspension of $\dl$ is 
in $\gor_a$.

In contrast to the case for Gorenstein* simplicial complexes there does
however not seem to be any formula for $\tvH_{d-1}(\dl)$ 
for $\dl$ in $\gor_a$
for instance in terms of $n$ and $d$.
\remfin

\section{Problems}

We pose the following two problems.

\begin{problem}
What are the possible $f$-vectors (or $h$-vectors) of the simplicial complexes
in the classes $\cler_a$ and $\cm_a$.
\end{problem}

This is likely to be a very difficult problem since any answer also would
include an answer to what the $h$-vectors of Gorenstein* simplicial complexes
are. However, any conjecture about this would be highly interesting since
it would contain as a subconjecture what the $h$-vectors of Gorenstein*
simplicial complexes are.

\begin{problem}
Construct simplicial complexes in the classes $\clerri_a$ and $\cmri_a$
for various parameters of $n,d,c$, and $a$.
\end{problem}

As has been pointed out this has been done in a number of particular
cases. When $a=0$ we have the bi-Cohen-Maculay simplicial complexes
constructed in \cite{FV}. When $a=1$ and $d=2c$ we have the cyclic polytopes
in $\cmri_1$, and when $a=c$, many Steiner systems $S(c,d,n)$ have been
constructed, which give simplicial complexes in $\clerri_a$.

\end{document}